\documentclass[12pt,letterpaper]{amsart}

\usepackage{graphicx}
\usepackage{times}
\usepackage{geometry}
\usepackage{fullpage}

\newtheorem{thm}{Theorem}

\newtheorem{Def}{Definition}
\newtheorem{prop}{Proposition}
\newtheorem{claim}{Claim}

\newcommand{\fin}{\mathcal{F}}
\newcommand{\bin}{\mathcal{B}}
\newcommand{\Rop}{\operatorname{Rop}}
\newcommand{\R}{\mathbb{R}}
\newcommand{\dcsd}{\operatorname{dcsd}}
\newcommand{\arcskip}{\operatorname{Skip}}
\newcommand{\Cr}{\operatorname{c}}

\title{Upper Bounds for Ropelength as a function of Crossing Number}

\author[Cantarella]{Jason Cantarella}
\address{Department of Mathematics, University of Georgia,
Athens, GA 30602}
\email{jason@math.uga.edu}

\author[Faber]{X.W. Faber}
\address{Department of Mathematics, University of Georgia,
Athens, GA 30602}
\email{faber@math.uga.edu}

\author[Mullikin]{Chad A. Mullikin}
\address{Department of Mathematics, University of Georgia,
Athens, GA 30602}
\email{chadm@math.uga.edu}

\begin{document}

\begin{abstract}
The paper provides bounds for the ropelength of a link in terms of the
crossing numbers of its prime components. As in earlier papers, the
bounds grow with the square of the crossing number; however, the
constant involved is a substantial improvement on previous
results. The proof depends essentially on writing links in terms of
their arc-presentations, and has as a key ingredient Bae and Park's
theorem that an $n$-crossing link has an arc-presentation with less
than or equal to $n+2$ arcs.
\end{abstract}

\keywords{ropelength, crossing number, arc-presentations, geometric knot theory}

\maketitle
\section{Introduction}

The {\em ropelength} of a space curve is defined to be the quotient of
its length by its {\em thickness}, where thickness is 
the radius of the largest embedded tubular neighborhood around the
curve. For a knot or link type $L$, we define the ropelength $\Rop(L)$ to be the
minimum ropelength of all curves with the given link type. This
minimum ropelength is a link invariant which measures the topological
complexity of the link, much like crossing number, or bridge number,
in classical knot theory.

It has been shown that every link type contains at least one $C^{1,1}$
{\em tight} representative which achieves this minimum ropelength
\cite{CKS,gmsvdm}. Much effort has been invested in the project of
finding lower bounds for the ropelength of various link types in terms
of classical topological invariants, such as the crossing number
\cite{buck,CKS,LSDR}. 

In this paper, we are interested in a converse problem: given a link
type $L$ of crossing number $\Cr(L)$, can we guarantee the existence of
a representative curve with ropelength less than some function of
$\Cr(L)$?  That is, can we find {\em upper} bounds on ropelength in terms of
crossing number? Our main theorem states the following:

\begin{thm}
If $L$ is a non-split link, then 
\begin{equation}
\Rop(L) \leq 1.64 \, \Cr(L)^2 + 7.69 \, \Cr(L) + 6.74.
\end{equation}
In particular, this bound holds for prime links.
\end{thm}

Our Theorem 2 gives similar bounds for composite links. 

Other groups (\cite{CKS,Joh}) have attacked this problem by
finding upper bounds on the number of edges required to embed a given
link $L$ in the unit lattice (the {\em lattice number} $k(L)$ of
the link), and then observing that $\Rop(L) < 2 k(L)$ \cite{DEJvR2}.
Both proofs rely on laying out a diagram of the knot as a graph in a
planar grid and then adding bridges to form overcrossings. In this
context, it has been observed that constructing a particular diagram
of a link with crossing number $\Cr(L)$ may require ropelength
$O(\Cr(L)^2)$ \cite{Joh}. These authors have obtained the weaker
bounds $\Rop(L) < 24 \Cr(L)^2$ \cite{CKS}, and $\Rop(L) < 25 \Cr(L)^2$
\cite{Joh}. Johnston's algorithm, like ours, produces an explicit
realization of the knot in space, while the approach of \cite{CKS} is
less constructive.

By contrast, our methods are more three-dimensional and are not based
on grid or lattice embeddings. Instead of using a planar diagram of a
knot, we base our construction on Peter Cromwell's idea of {\em
arc-presentations} \cite{Cr}. It is curious that our methods, too, seem to be
essentially of order $\Cr(L)^2$.  While we believe that bounds with a
slower order of growth must be attainable, it is becoming clear that the
problem of constructing such bounds is likely to be challenging.

\section{The Definition of Ropelength}

The ropelength of a curve is defined to be the quotient of length by
the radius of the largest embedded tubular neighborhood around the
curve. This radius is called the {\em thickness} of the
curve. For~$C^2$ curves, this radius is locally controlled by
curvature and globally controlled by distances of self-approach
between various regions of the curve. Formally, we write

\begin{Def} \label{def:thickness}
The {\em thickness} of a $C^2$ curve $c$ is given by
\begin{equation}
\tau[c] := \min \left\{ \min_s \frac{1}{\kappa(s)}, \frac{\dcsd(c)}{2} \right\},
\end{equation}
where $\kappa(s)$ is the curvature of $c$ at $s$, and $\dcsd(c)$ is
the shortest {\em doubly-critical self-distance} of $c$; that is, the
length of the shortest chord of $c$ which is perpendicular to the
tangent vector $c'$ at both endpoints.
\end{Def}

We can extend this defintion to $C^{1,1}$ curves by adjusting our idea
of the radius of curvature as follows (c.f. \cite{CKS}):

\begin{Def} \label{def:fixit}
Let $s$ be a point on a $C^{1,1}$ curve. Consider a decreasing
sequence of open neighborhoods $U_n$ of $s$. The {\em infimal radius
of curvature} at $s$ is given by
\begin{equation}
  \inf_{U_n} \left\{ \inf_{t \in U_n} \frac{1}{\kappa(t)}\right\},
\end{equation}
where the inner infimum is restricted to $t$ in $U_n$ such that
$\kappa(t)$ exists.
\end{Def}

\begin{figure}[ht!]
\centerline{\includegraphics[height=1.4in]{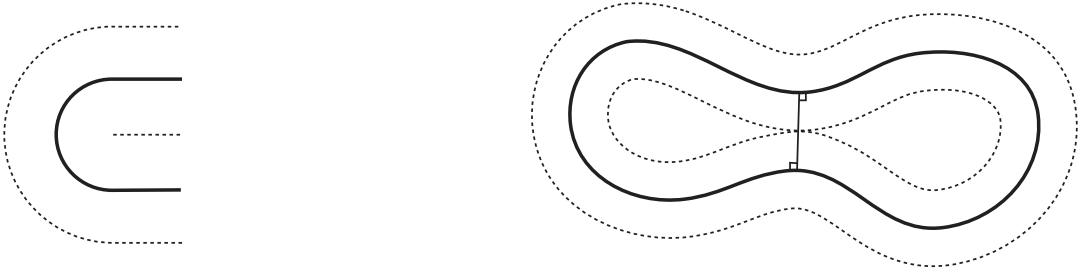}}
\caption{These are two curves of unit thickness in the plane with their
largest embedded tubular neighborhoods. In the left curve, thickness
is controlled by curvature while in the right curve, thickness is
controlled by the length of the doubly-critical chord shown.}
\label{fig:thick}
\end{figure}

Figure~\ref{fig:thick} shows examples of curves where thickness is controlled by
curvature and by the doubly-critical self-distance.

Gonzalez and Maddocks have given another definition of thickness which
looks somewhat less natural, but is often more useful. (See \cite{GM}
for details). Another useful way to look at thickness comes from
Federer's notion of {\em reach}, which agrees with the thickness for
curves \cite{federer}.

\begin{Def} 
The {\em reach} of a set $S$ inside $\R^n$ is the greatest non-negative $r$
so that each point within distance $r$ of $S$ has a unique nearest neighbor
in~$S$.
\end{Def} 

\section{Arc-presentations}

We start with a definition:

\begin{Def} 
\label{def:arcpres}
An {\em arc-presentation} of a link $L$ is an embedding of $L$ in a
finite collection of $\alpha$ open half-planes arrayed around a common
axis, or binding, so that the intersection of $L$ with each half-plane
is a single simple arc. The number of half-planes $\alpha$ is called
the {\em arc-index} of the arc-presentation.  The minimal arc-index
over all arc-presentations of a link $L$ is an invariant of the link type.
\end{Def}

By isotopy, we can arrange that $L$ intersects the axis only at the
points $1, \dots, \alpha$. We call these the {\em levels} of the
arc-presentation.  Such an arc-presentation is then specified by
combinatorial data: a collection of $\alpha$ triples in the form
$(x_i,y_i,\theta_i)$, where each denotes an arc from level $x_i$ to
level $y_i$ on the half-plane at angle $\theta_i$ around the axis.

\begin{center}
\begin{figure}[ht!]
\begin{minipage}[b][3in][t]{5cm}
\includegraphics[height=3in]{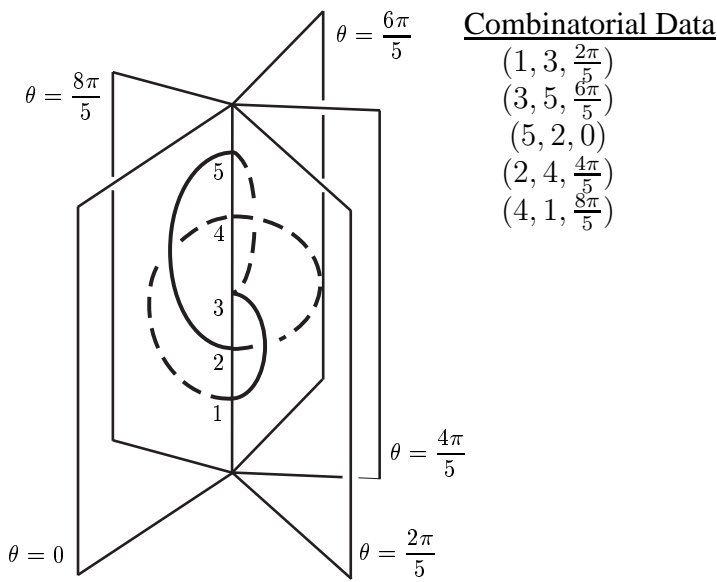}
\end{minipage}
\begin{minipage}[b][3in][t]{4cm}
\begin{tabular}[b]{ c }
\qquad \underline{Combinatorial Data} \\
$(1, 3, \frac{2\pi}{5})$ \\
$(3, 5, \frac{6\pi}{5})$ \\
$(5, 2, 0)$ \\
$(2, 4, \frac{4\pi}{5})$ \\
$(4, 1, \frac{8\pi}{5})$ 
\end{tabular}
\end{minipage}
\caption{This figure shows an arc-presentation for a trefoil
knot. The presentation has arc-index $5$. To the right we see
the combinatorial data which describes this arc-presentation: $5$
triples in the form $(x_i,y_i,\theta_i)$, each indicating an arc from
level $x_i$ to level $y_i$ on page $\theta_i$ of the ``$5$-page book''
shown on the left.}
\label{ArcPresTref}
\end{figure}
\end{center}

Figure~\ref{ArcPresTref} shows an arc-presentation for the trefoil and the
corresponding set of triples.

We will assemble our ropelength bounds from two ingredients. First, we
define a notion of the total distance travelled by the arcs in an
arc-presentation:

\begin{Def}
\label{def:skip}
The {\em total skip} of an arc-presentation $A$, denoted $\arcskip(A)$, is
\begin{equation}
\arcskip(A) = \sum_{i=1}^{\alpha} |x_i - y_i|.
\end{equation}
\end{Def}

For a given arc-presentation we can construct a realization of the
knot in space with ropelength bounded in terms of~$\arcskip(A)$
and~$\alpha$:

\begin{prop} 
An arc-presentation $A$ composed of $\alpha$ half-planes can be
realized with ropelength smaller than
\begin{equation}
\frac{2\alpha}{\tan (\pi / \alpha)} + (\pi - 2) \alpha + 2 \arcskip(A).
\end{equation}
\label{prop:skipbound}
\end{prop}

\begin{figure}[ht!]
    \begin{tabular}{ccc}
    \includegraphics[height=2.5in]{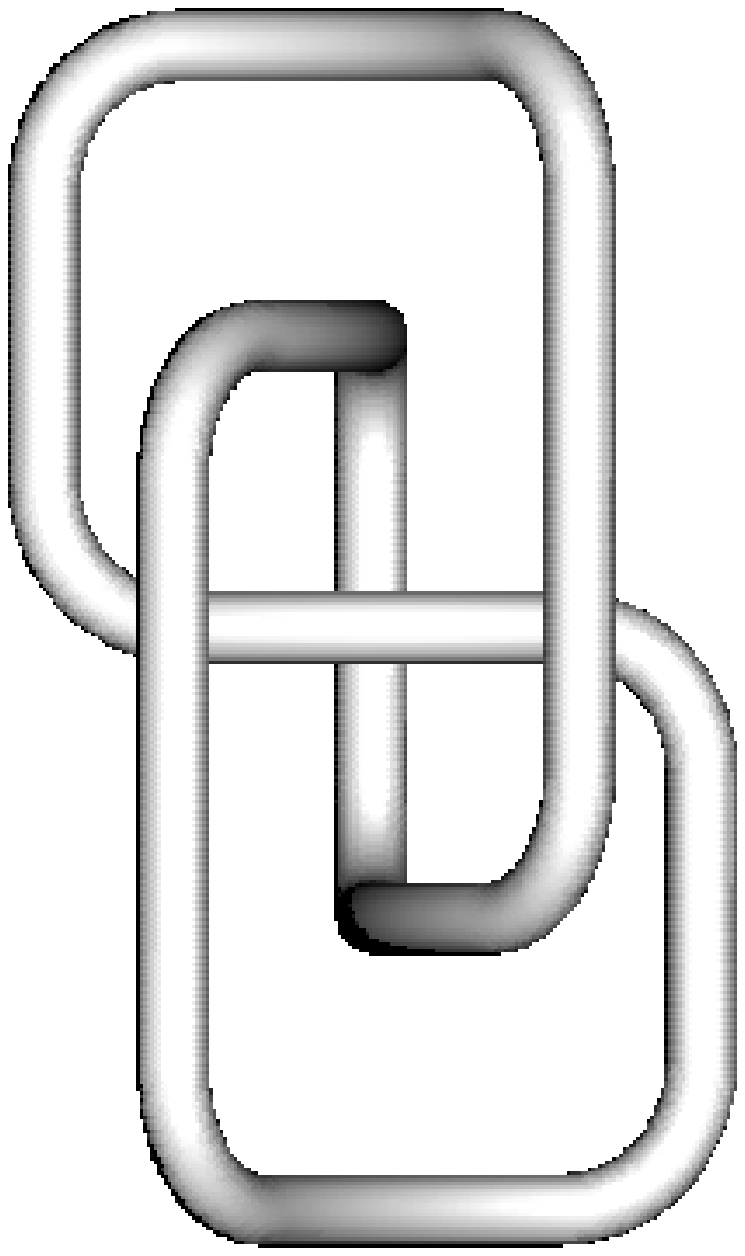} & 
    \hspace{1in} &
    \includegraphics[height=2.5in]{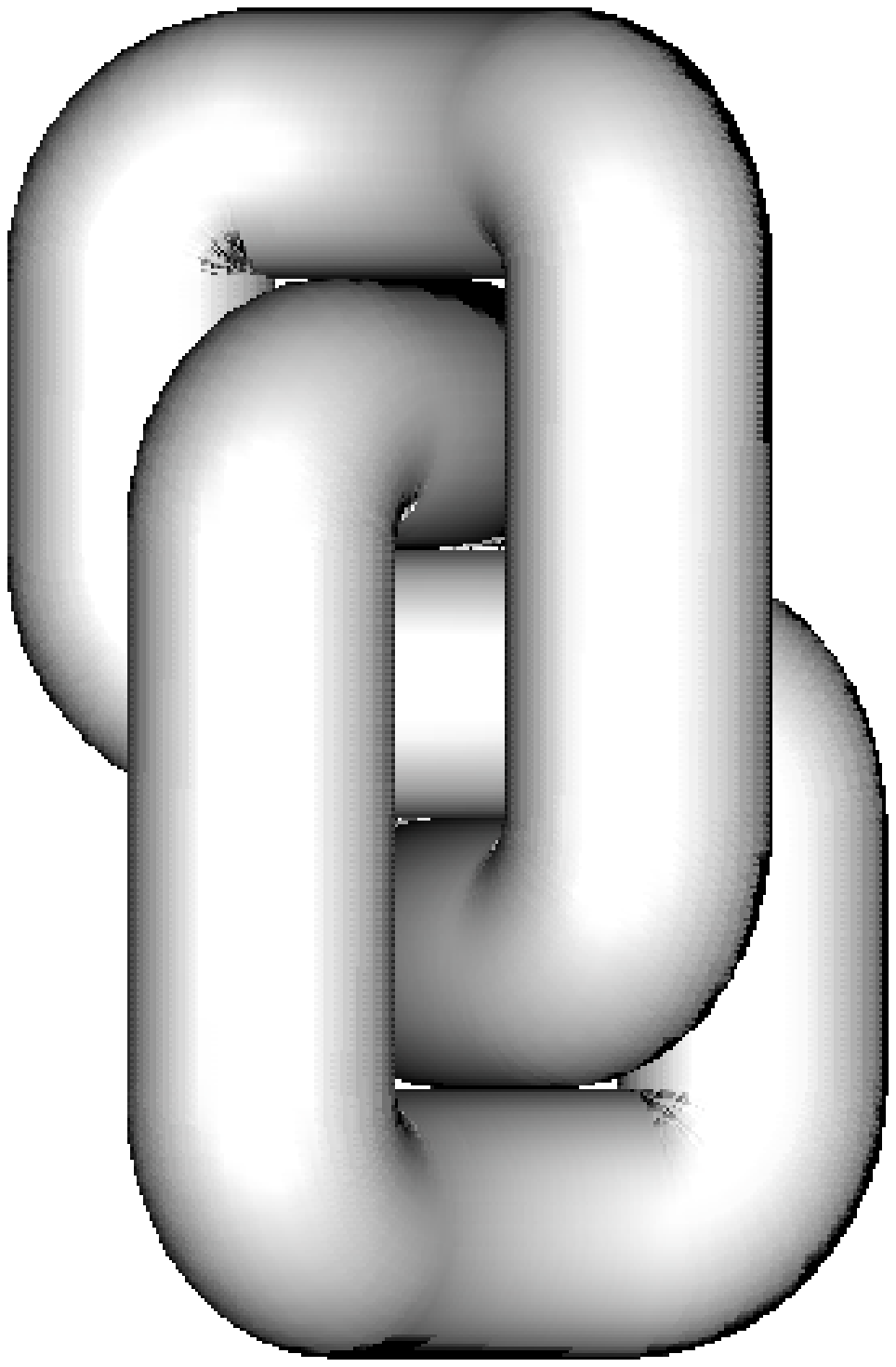} 
    \\
    \includegraphics[height=2.5in]{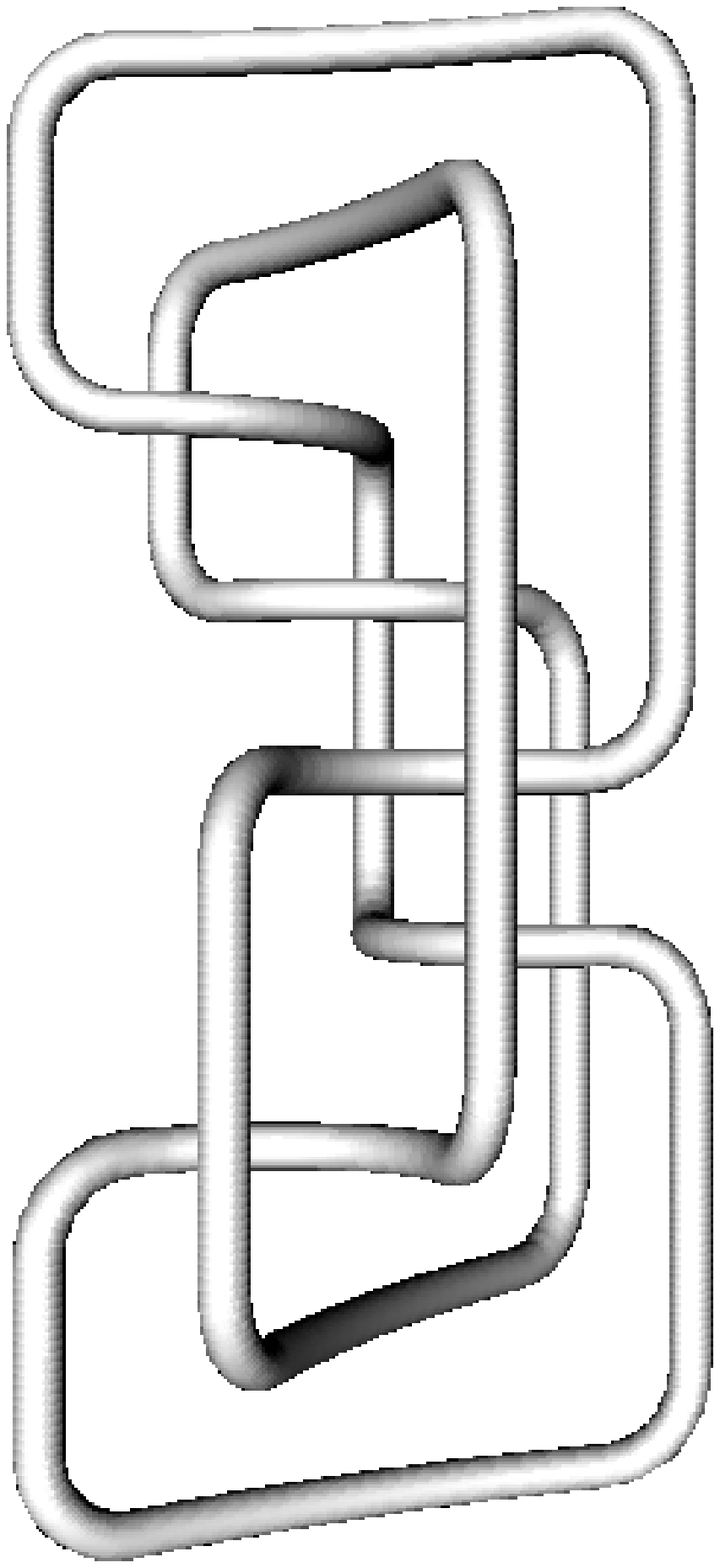} &
    \hspace{1in} &
    \includegraphics[height=2.5in]{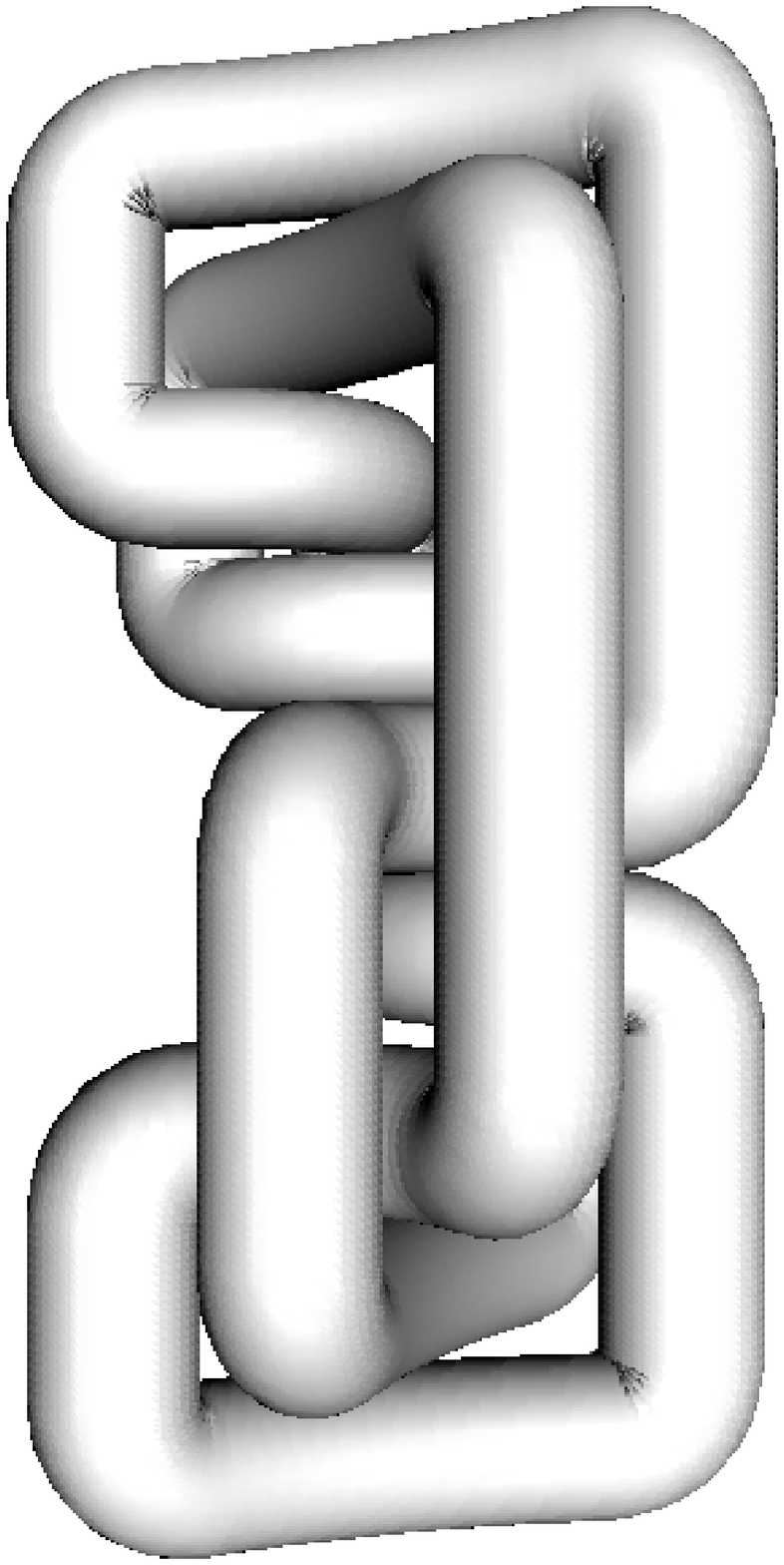}\\
    \end{tabular}
\caption{Here we see a trefoil knot (top left) and a $7_1$ knot
(bottom left) together with the tubular neighborhoods around them
constructed by Proposition~1.  Our trefoil knot
appears much tighter: its ropelength~($43.47$) is
proportionally closer to the minimum ropelength for its knot
type~($32.66$) than the ropelength of our $7_1$~knot~($97.05$) is to
the minimum for its knot type~($61.40$~\cite{sdkp}).}
\label{fig:examples}
\end{figure}

For the arc-presentation of the trefoil in Figure~\ref{ArcPresTref},
we have $\alpha = 5$ and $\arcskip(A) = 12$; so
Proposition~\ref{prop:skipbound} yields an upper bound on the
ropelength of the trefoil of about~$43.47$. Numerical experiments
estimate the ropelength of the tight trefoil to be about
$32.66$~\cite{sdkp}, so the slack in our estimate is 
about $33\% $ of the total value. Figure~\ref{fig:examples} shows the
tubular neighborhoods of this trefoil knot and an arc-presentation of
the knot~$7_1$ as realized by the algorithm of Proposition~1.

Further, if we can bound $\arcskip(A)$ in general, we will be able
to draw conclusions about the ropelength of an arbitrary link. A
combinatorial argument yields:

\begin{prop} If an arc-presentation $A$ has arc-index $\alpha$, then
\begin{equation}
\arcskip(A) \leq\begin{cases} \frac{\alpha^2 - 1}{2}& \text{if
$\alpha$ is odd,} \\ \frac{\alpha^2}{2}& \text{if $\alpha$ is even.}
\end{cases}
\end{equation}
This bound is sharp.
\end{prop}

It is shown in \cite{baepark} that any non-split link $L$ admits an
arc-presentation with $\alpha \le \Cr(L) + 2$.  This result, when coupled
with the previous two propositions, gives Theorem 1.  We obtain an
even stronger statement for composite links:

\begin{thm} 
If $L$ is a non-split composite link with prime components $L_1, L_2, \ldots,
L_n$, then
\begin{equation}
\Rop(L) \leq 1.64\sum_{i=1}^n \Cr(L_i)^2 + 7.69\sum_{i=1}^n \Cr(L_i) +
6.74n.
\end{equation}
\end{thm}

\section{Proofs of the key propositions and theorems}

\begin{proof}[Proof of Proposition 1]
We would like to take an arc-presentation $A$ for $L$ as a template for 
constructing an embedding of $L$ with unit
thickness.  We will then bound the length of this embedding in terms
of the arc-index and the total skip of $A$.

We begin by constructing a right regular polygonal prism $P \times
[0,2\alpha]$, where $P$ is a regular polygon with $\alpha$ sides of
length $2$. This prism will serve as the binding of $A$; each
vertical face of the prism will correspond to an open half-plane in the
arc-presentation $A$. We divide the prism vertically into $\alpha$
{\em floors}, each a prism of height $2$, which will represent the
$\alpha$ {\em levels} of the arc-presentation $A$.

We can now construct a link isotopic to $L$. First, represent the arcs
of $A$ by $\alpha$ handles outside the prism which join different
floors on the same vertical face. We will refer to these handles as
{\em fins}. Next, add $\alpha$ circular sections inside the prism
which join different vertical faces on the same floor. These sections
represent the junctions between arcs on the binding of the open book
described by~$A$.

We must show that this construction can be accomplished with a unit
thickness curve and then compute the length of that curve.

\subsection{The Fins}

Let us denote the fins $\fin_1, \dots, \fin_\alpha$.  Each fin
consists of two quarter-circles of unit radius, joined by a straight
vertical segment. Each fin joins two points on a vertical face of the
prism and is contained in a $2\alpha \times 2 \times 2$ rectangular
box extending radially from a vertical face of the prism.

Since the $\fin_i$'s stay outside the prism and each is contained in a
different box, the tubes around the fins are pairwise disjoint, and
disjoint from the tubes surrounding regions of the curve inside the
prism. Given that each fin has curvature bounded above by $1$ and no
doubly-critical chords, this means that the fins can be constructed
with a unit-thickness curve.

\begin{claim}
If $\Rop(\fin_i)$ denotes the length of the segment of the
curve~$\fin_i$, then
\begin{equation}
\sum_{i=1}^{\alpha} \Rop(\fin_i) = (\pi - 2) \alpha + 2\,\arcskip(A).
\label{eq:rfin}
\end{equation}
\end{claim}

\begin{proof} Suppose that $\fin_i$ travels from floor $x_i$ to floor $y_i$ of the
prism. The total vertical distance covered by the fin is $2|x_i -
y_i|$ (recall that each floor has height $2$). However, the
quarter-circles on each end of the fin cover a vertical distance of
$2$ units. Thus, the straight segment has length $2|x_i~-~y_i| - 2$,
and the total length of the fin is $\pi - 2 + 2 |x_i - y_i|$. Summing
over $i = 1, \dots, \alpha$ and using Definition~\ref{def:skip}
proves the claim. 
\end{proof}

\subsection{The Binding Prism}

We denote the sections of the curve inside each floor of the binding
prism by $\bin_1, \dots, \bin_\alpha$. Each $\bin_i$ is a circular
arc joining the midpoints of two edges of the regular polygon which 
is the cross-section of the prism as shown in Figure~\ref{fig:paths}.

\begin{figure}[ht!]
\includegraphics[width=1.5in]{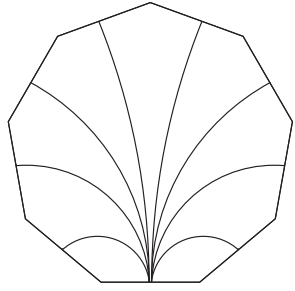}
\caption{The sections of our curve $\bin_i$ within the binding prism
are circular arcs joining the midpoints of edges of the cross-section
of the prism. The plane of this picture is
located in the center of a floor of the prism.}
\label{fig:paths}
\end{figure}

Because the sides of the polygon have length~$2$, each of these is an
arc of a circle of radius at least one; so each arc has curvature bounded
above by one. Further, since each floor has height~$2$ and only one
$\bin_i$ lies in each floor, the tubes around each of the~$\bin_i$ 
are disjoint. Thus these~$\bin_i$ can be constructed with a tube 
of unit thickness. 

\begin{claim}
If $\Rop(\bin_i)$ denotes the length of the segment of the curve~$\bin_i$,
then 
\begin{equation}
\sum_{i=1}^\alpha \Rop(\bin_i) \leq \frac{2\alpha}{\tan(\pi/\alpha)}.
\label{eq:rbin}
\end{equation}
\end{claim}

\begin{proof}
Each of these circular arcs is contained in a sector of the circle
inscribed within the polygonal cross-section of the prism as shown in
Figure~\ref{fig:bounding}.  Since each arc is convex, its length
is bounded above by the diameter of the inscribed circle.  This diameter is
exactly~$2\cot(\pi/\alpha)$. Summing over $i = 1, \dots, \alpha$
proves the claim.
\end{proof}

\begin{figure}[ht!]
\centerline{\includegraphics[width=1.5in]{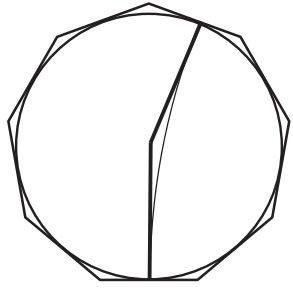}}
\caption{Each of the paths $\bin_i$ through the floors of the
binding prism is a circular arc connecting two sides of the polygon which
is that prism's cross-section. Here we see that each of these arcs is
contained within a sector of the circle inscribed within that
polygon. Since each arc is a convex curve, this means that its length
is bounded by the length of the two radii which bound the sector. That is, it is bounded
by the diameter of the inscribed circle.}
\label{fig:bounding}
\end{figure}
 
Combining Claims~1 and 2 yields the statement of Proposition~1.
\end{proof}

\begin{proof}[Proof of Proposition 2]
Our job is to find an upper bound for $\arcskip(A) = \sum_{i=1}^\alpha
|x_i - y_i|.$ We first observe that the difference $|x_i - y_i|$ is
one unit larger than the number of levels skipped over. For example,
jumping from level~$3$ to level~$6$, a difference of $3$ levels, skips
the fourth and fifth levels. Thus, we can rewrite the sum
\begin{equation}
\arcskip(A) = \alpha + \sum_{i=1}^{\alpha} \{\text{number of levels skipped
by the arc $(x_i,y_i,\theta_i)$}\}.  
\end{equation}
Notice that any level~$j$ contributes to the above sum exactly when it is skipped over.
We can rewrite our sum in terms of $j$ as
\begin{equation}
\begin{aligned}
\arcskip(A) &= \alpha + \sum_{j=1}^{\alpha} 
               \{\text{number of times level $j$ is skipped}\} \\
            &= \alpha + \sum_{j=1}^{\lfloor \alpha / 2 \rfloor} 
               \{\text{number of times level $j$ is skipped}\} \\
            & \qquad  + \sum_{j=0}^{\alpha -\lfloor \alpha /2 \rfloor-1}
               \{\text{number of times level $\alpha - j$ is skipped}\},
\end{aligned}
\label{eq:skips}
\end{equation}
where in the final equality we have split the second half of the sum
off and let $j\mapsto \alpha -j$.

Now we bound the number of times level $j$ is skipped over.  The only
way to hop over $j$ from a higher level is to land on a lower
level. There are $j-1$ levels below the $j$th on which such a jump can
land.  Further, each of these levels can act as a launch pad for a
jump back up which crosses the $j$th level again.  This gives at most
$2(j-1)$ skips over level $j$.  Similarly, the number of times we can
skip over the $\alpha - j$th level is twice the number of levels above
it, or~$2j$.

For even $\alpha$, these estimates are sharp (as we will see below).
However, when level $\alpha-j$ is the central level of an arc-presentation
with~$2k+1$ levels ($j=k=\frac{\alpha-1}{2}$), the situation is slightly different. Here all of
the~$j$ levels above the middle cannot be initial and terminal
levels of arcs which skip level~$\alpha-j$. For if so, then no arcs land on
level $\alpha-j$, and we could have eliminated level $\alpha-j$ from the original
arc-presentation. Thus level $\alpha-j$ is skipped at most $2j-1=\alpha-2$ times.

Inserting these bounds into Equation~\ref{eq:skips}, we apply
the sum formulae for arithmatic progressions. When $\alpha$ is odd,
we get
\begin{equation}
\arcskip(A) \le \alpha + \sum_{j=1}^{\frac{\alpha-1}{2}} 2(j-1)
+ \sum_{j=0}^{\frac{\alpha -3}{2}} 2j + (\alpha - 2)
=\frac{\alpha^2 - 1}{2}. 
\end{equation}
If $\alpha$ is even, the proof is similar.

We now construct arc-presentations which show that these results are
sharp.  Consider the arc-presentation with even arc-index $\alpha=2k$
described by the data 
\begin{equation*}
\begin{aligned}
(\alpha,\alpha/2,\theta_1), (\alpha/2, \alpha-&1,\theta_2),(\alpha-1,\alpha/2-1,\theta_3), 
(\alpha/2-1, \alpha-2, \theta_4) , \\ 
& \ldots, (\alpha/2+1,1,\theta_{2k-1}), (1,\alpha,\theta_{2k}).
\end{aligned}
\end{equation*}
If we add up the lengths of the jumps, we get 
\begin{equation}
\begin{aligned}
\arcskip(A) = \alpha^2 / 2.
\end{aligned}
\end{equation}
The same approach yields a realization of $A$ so that $\arcskip(A) = \frac{\alpha^2 - 1}{2}$ for odd $\alpha$.

\end{proof}

\begin{proof}[Proof of Theorem 1]
Taylor's theorem gives the approximation $\frac{1}{\tan(x)} \leq 1/x -
x/3$ for $x>0$.  Via Propositions 1 and 2 we gather that
\begin{equation}
\begin{aligned}
\Rop(L) &\leq \frac{2\alpha}{\tan(\pi / \alpha)} + (\pi - 2)\alpha + \alpha^2 \\
&\leq (2/{\pi} + 1)\alpha^2 + (\pi - 2)\alpha - 2\pi/3.
\end{aligned}
\end{equation}
By Bae and Park \cite{baepark}, for any non-split link $L$ there exists an
arc-presentation with $\alpha \leq \Cr(L) + 2$.  Inserting this into
the above bound for ropelength yields
\begin{equation}
\Rop(L)\leq (2/\pi +1) \Cr(L)^2 + (8/\pi + 2 + \pi) \Cr(L) +(8/\pi + 4\pi /3),
\end{equation}
and each of these constants evaluates to something smaller than the
approximations given in the statement of the theorem.  To gain the
final remark in the theorem, we note that any prime link $L$ is
non-split (otherwise it would consist of split compnents $L_1$ and
$L_2$ and would admit the nontrivial factors $L_1$ and $L_2$ union a
split unknot).
\end{proof}



\begin{proof}[Proof of Theorem 2]
The strategy for this proof is to arrange the prime components of our
composite link so that we can make use of the bounds given by
Theorem~1. So suppose that we have found arc-presentations with
minimal arc-index for these components and embedded them as
unit-thickness curves $L_1, \dots, L_n$ according to the algorithm of
Proposition~1.

We will now prove that for any links $L_1$ and $L_2$, constructed by the
algorithm of Proposition~1, we can construct a curve $L_1 \# L_2$ with
ropelength less than or equal to $\Rop(L_1) + \Rop(L_2)$.  This is all
that is required to complete the proof of our Theorem since the bound
in the statement is just the sum of the bounds obtained for the $L_i$
by Theorem~1.

We begin by preparing $L_1$ and $L_2$. The top floor of $L_1$ contains
only a single horizontal circular arc joining the centers of two sides of the binding
prism. Since no fins jump over this level, we may rotate
these quarter-circles to face one another and replace the horizontal circular
arc with a horizontal line segment of shorter length without changing thickness
or knot type. We do the same for the bottom floor of $L_2$. This
procedure is shown in Figure~\ref{fig:snorkel_cs}.

\begin{figure}[ht!]
\begin{center}
\includegraphics[height=1in]{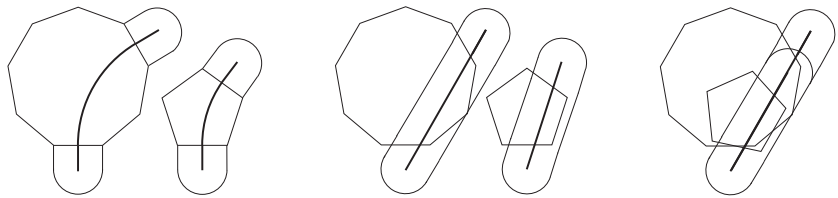}
\end{center}
\caption{We look down on a knot of arc-index 5, whose binding prism is
  shown by the small pentagon at right, preparing to be joined to a
  knot of arc-index 9, whose binding prism is shown by the large
  nonagon at left. The leftmost pair of figures shows the original
  position of the top and bottom arcs of these knots, while the middle
  pair of figures shows these arcs ``straightened'' to prepare for the
  connect sum.  The rightmost pair of figures shows the two binding
  prisms in the correct relative position for the connect sum.}
\label{fig:snorkel_cs}
\end{figure}

We now arrange $L_1$ and $L_2$ in space so that the horizontal
segments are colinear and share an endpoint. If we keep each oriented
so that its floors are horizontal, the only overlap between the tubes
surrounding each curve occurs on the shared floor. At the shared
endpoint, we may delete two quarter-circles and replace them with a
vertical line segment of length $2$. We could
keep track of this savings and get a slightly better constant term in
the statement of Theorem~2. For each prime component we add, we save
$\pi - 2$ in length.

Handling the other endpoints of the curve will prove to be a little
more work.  We may assume that both line segments lie along the
x-axis with the shared endpoint at the origin. Suppose $L_2$'s
segment has length $\ell_2$, while $L_1$'s segment has the smaller
length~$\ell_1$.

We now rotate the remaining vertical quarter-circle of $L_1$ to face
the corresponding quarter-circle of $L_2$. If $\ell_1 \leq \ell_2 -
2$, we may replace both horizontal line segments with a single,
shorter horizontal line segment joining the ends of these vertical
quarter circles to obtain the desired curve. See Figure~\ref{fig:csone}.

\begin{figure}[ht]
\centerline{\includegraphics[height=0.75in]{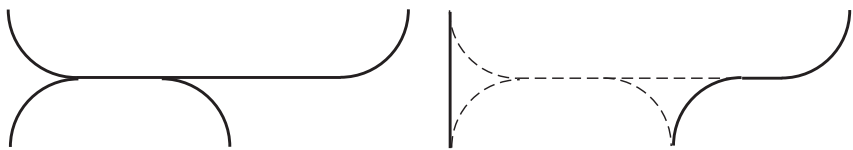}}
\caption{This figure shows the extreme arcs of the two components of
  the connect sum, straightened, and aligned with one another on the
  left. On the right, we see the new curve. Two quarter-circles on the
  left have been replaced with a straight line segment; the lower
  quarter-circle has been rotated to face right; the lower horizontal
  segment (of length $\ell_1$) has been deleted; and the upper
  horizontal segment (of length $\ell_2$) has been replaced by a
  horizontal segment of length $\ell_2 - \ell_1 - 2$. Since these
  changes all reduce length, the curve on the right is strictly
  shorter.}
\label{fig:csone}
\end{figure}

If $\ell_1 > \ell_2 - 2$, we cannot simply connect the endpoints of
the quarter-circles after rotating the lower quarter-circle to face
right. The resulting curve would have cusps on both ends. We solve
this problem by finding a line tangent to both circles and following
the composite path shown in Figure~\ref{fig:csthree}.

\begin{figure}
\centerline{\includegraphics[height=0.75in]{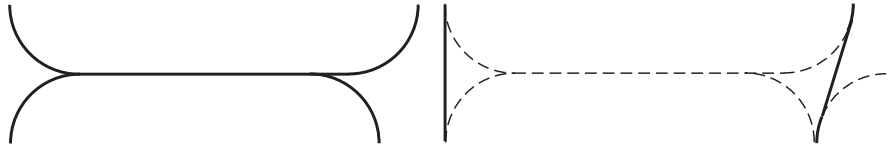}}
\caption{This figure shows the two extreme arcs of the components in 
  the case where $\ell_1 > \ell_2-2$. When we
  rotate the lower quarter-circle to face right, it cannot be joined
  by a horizontal straight straight to the upper quarter-circle to
  create a $C^{1,1}$
  curve; instead we find the diagonal line tangent to both
  quarter-circles and follow the composite path shown.}
\label{fig:csthree}
\end{figure}

It is less obvious that these changes reduce length. To see that they
do, we consider the diagonal line tangent to both circles shown in
Figure~\ref{fig:csthree}. Since both circles are also tangent to a
horizontal line, by symmetry this horizontal line cuts the diagonal
line in half. Consider Figure~\ref{fig:csfour}. We need only show that
half of the diagonal line (labelled $x$ in the Figure) is shorter than
the portion of the quarter-circle it replaces (twice the angle
$\theta$).

\begin{figure}
\centerline{\includegraphics[height=1.5in]{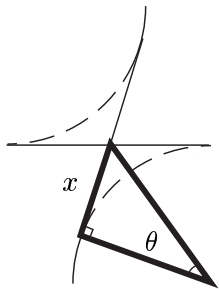}}
\caption{This detailed Figure enlarges the right-hand side of
  Figure~\ref{fig:csthree}. Consider the triangle with the following
  vertices: the point of tangency of the diagonal segment with the lower
  circle, the center of the lower circle, and the midpoint of the
  diagonal segment.
. The portion of the lower quarter-circle replaced by this half
  of the line segment has length $2\theta$ (again by symmetry). The length of
  this portion of the line segment is given by $x$.}
\label{fig:csfour}
\end{figure}

Since the lower quarter-circle has unit radius, this amounts to
proving that $\tan\theta \leq 2\theta$ for $0 \leq \theta \leq \pi/4$.
This is shown by a simple computation.

Since the resulting curve remains $C^{1,1}$, is still of unit thickness, and
has less length than the total length of the initial curves, this
completes the proof.
\end{proof}

\begin{figure}[ht!]
\begin{center}
\begin{tabular}{ccc}
\includegraphics[height=2.3in]{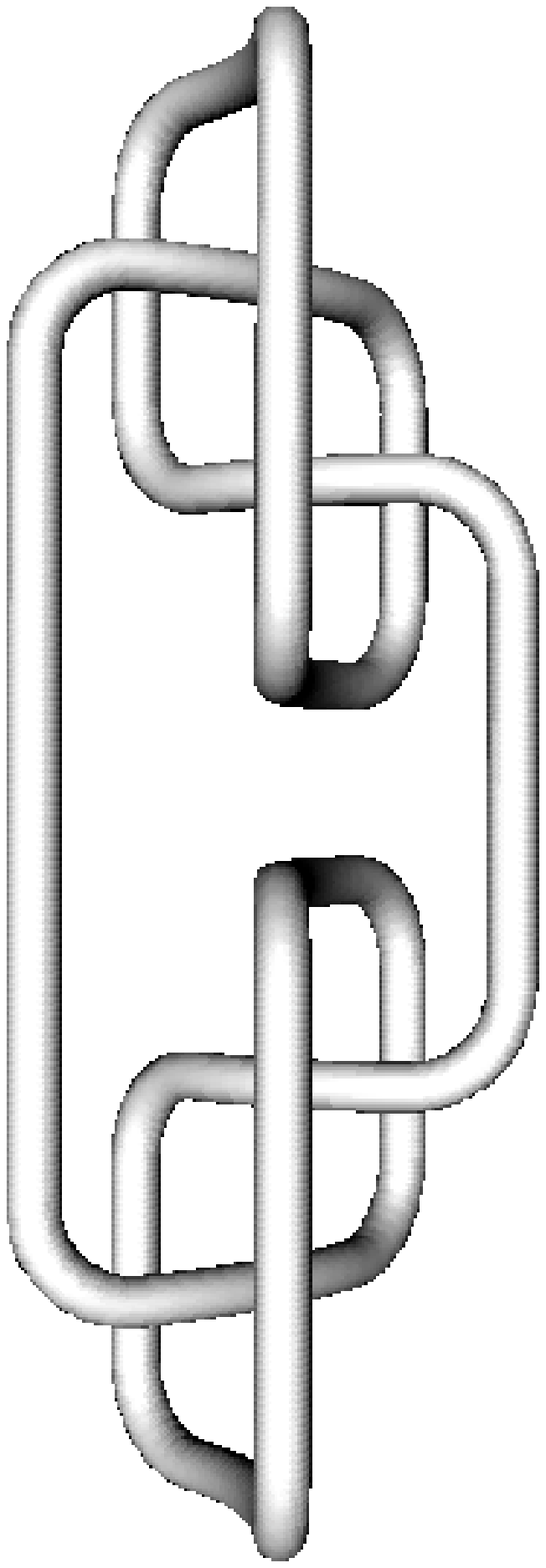}
&
\hspace{1in}
&
\includegraphics[height=2.3in]{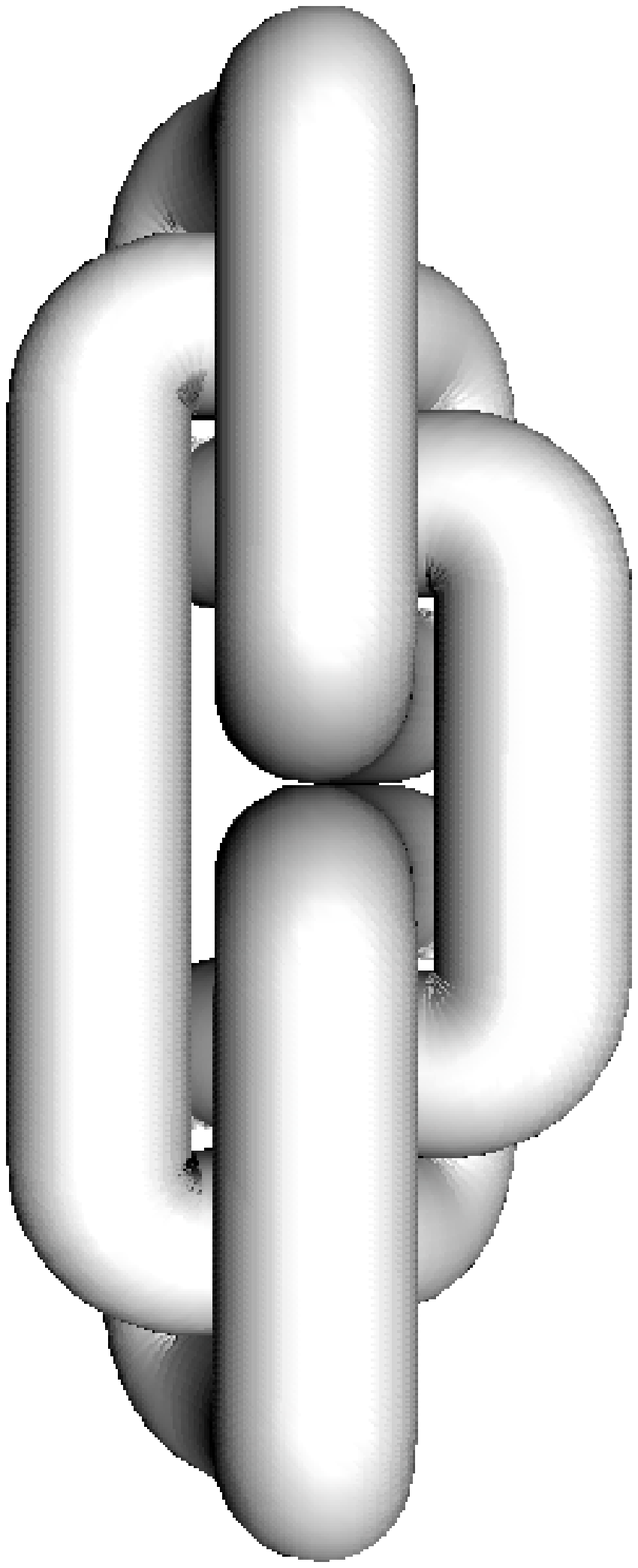}
\end{tabular}
\end{center}

\caption{Here we see the results of the construction of Theorem~2. Two
  mirror-image trefoil knots, generated by the method of Proposition~1
  from the arc-presentation given in Figure~\ref{ArcPresTref}, have
  been joined by the methods of Theorem~2 to obtain the composite knot
  $3_1 \# \overline{3_1}$.}

\label{fig:ready_to_go}
\end{figure}

An example of this construction is shown in
Figure~\ref{fig:ready_to_go}.

\section{Acknowledgements}

We would like to thank many of our colleagues for helpful
conversations on these topics over several years, including Yuanan
Diao, Heather Johnston, Rob Kusner and John Sullivan. We are indebted 
to the VIGRE program at the University of Georgia,
and to the other members of our VIGRE research group: Ted Ashton,
Kenny Little, Heunggi Park,
Darren Wolford, and Nancy Wrinkle. In addition, Cantarella would
like to acknowledge the generous support of the National Science 
Foundation Postdoctoral Research Fellowship Program under grant number
DMS-99-02397, and Faber and Mullikin would like to recognize
the gracious support of the VIGRE fellowship program at UGA.

\bibliographystyle{plain}
\bibliography{thick}

\begin{thebibliography}{10}

\bibitem{baepark}
Yongju Bae and Chan-Young Park.
\newblock An upper bound of arc index of links.
\newblock {\em Math. Proc. Cambridge Philos. Soc.}, 129(3):491--500, 2000.

\bibitem{buck}
Greg Buck and Jon Simon.
\newblock Thickness and crossing number of knots.
\newblock {\em Topol. Appl.}, 91(3):245--257, 1999.

\bibitem{CKS}
Jason Cantarella, Robert Kusner, and John Sullivan.
\newblock On the minimum ropelength of knots and links.
\newblock {\em Inventiones Mathematicae}, 150(2):257--286, 2002.

\bibitem{Cr}
Peter~R. Cromwell.
\newblock Arc presentations of knots and links.
\newblock In {\em Knot theory (Warsaw, 1995)}, pages 57--64. Polish Acad. Sci.,
  Warsaw, 1998.

\bibitem{DEJvR2}
Yuanan Diao, Claus Ernst, and E.J.~Janse van Rensburg.
\newblock Upper bounds on linking number of thick links.
\newblock Preprint, 2002.

\bibitem{federer}
Herbert Federer.
\newblock Curvature measures.
\newblock {\em Trans. Amer. Math. Soc.}, 93:418--491, 1959.

\bibitem{gmsvdm}
O.~Gonzalez, J.~H. Maddocks, F.~Schuricht, and H.~von~der Mosel.
\newblock Global curvature and self-contact of nonlinearly elastic curves and
  rods.
\newblock {\em Calc. Var. Partial Differential Equations}, 14(1):29--68, 2002.

\bibitem{GM}
Oscar Gonzalez and John~H. Maddocks.
\newblock Global curvature, thickness, and the ideal shapes of knots.
\newblock {\em Proc. Nat. Acad. Sci. (USA)}, 96:4769--4773, 1999.

\bibitem{Joh}
Heather Johnston.
\newblock An upper bound on the minimal edge number of an n-crossing lattice
  knot.
\newblock Preprint.

\bibitem{LSDR}
Richard~A. Litherland, Jon Simon, Oguz Durumeric, and Eric Rawdon.
\newblock Thickness of knots.
\newblock {\em Topol. Appl.}, 91(3):233--244, 1999.

\bibitem{sdkp}
Andrzej Stasiak, Jacques Dubochet, Vsevolod Katritch, and Piotr Pieranski.
\newblock Ideal knots and their relation to the physics of real knots.
\newblock In {\em Ideal knots}, volume~19 of {\em Ser. Knots Everything}, pages
  1--19. World Sci. Publishing, River Edge, NJ, 1998.

\end{thebibliography}

\end{document}